# On Pre-γ-*I*-Open Sets In Ideal Topological Spaces

Hariwan Zikri Ibrahim
Department of Mathematics, Faculty of Science, University of Zakho, Kurdistan Region-Iraq

**ABSTRACT**
In this paper, the author introduce and study the notion of pre-γ-*I*-open sets in ideal topological space.

## 1. INTRODUCTION

In 1992, Jankovic and Hamlett introduced the notion of *I*-open sets in topological spaces via ideals. Dontchevin 1999 introduced pre-*I*-open sets, Kasaharain 1979 defined an operation α on a topological space to introduce α-closed graphs. Following the same technique, Ogata in 1991defined an operation γ on a topological space and introduced γ-open sets. In this paper, some relationships of pre-γ-*I*-open, pre-*I*-open, preopen, pre-γ-open, γ-p-open, γ-preopen, *I*-open, $\delta_I$-open, R-*I*-open, α-*I*-open, semi-*I*-open, b-*I*-open and weakly *I*-local closed sets in ideal topological spaces are discussed.

## 2. PRELIMINARIES

Throughout this paper, $(X, \tau)$ and $(Y, \sigma)$ stand for topological spaces with no separation axioms assumed unless otherwise stated. For a subset $A$ of $X$, the closure of $A$ and the interior of $A$ will be denoted by $Cl(A)$ and $Int(A)$, respectively. Let $(X, \tau)$ be a topological space and $A$ a subset of $X$. A subset $A$ of a space $(X, \tau)$ is said to be regular open [N. V. Velicko, 1968] if $A = Int(Cl(A))$. $A$ is called δ-open [N. V. Velicko, 1968] if for each $x \in A$ there exists a regular open set $G$ such that $x \in G \subseteq A$. An operation γ [S. Kasahara, 1979] on a topology τ is a mapping from τ in to power set $P(X)$ of $X$ such that $V \subseteq \gamma(V)$ for each $V \in \tau$, where γ(V) denotes the value of γ at $V$. A subset $A$ of $X$ with an operation γ on τ is called γ-open [H. Ogata, 1991] if for each $x \in A$, there exists an open set $U$ such that $x \in U$ and $\gamma(U) \subseteq A$. Then, $\tau_\gamma$ denotes the set of all γ-open set in $X$. Clearly $\tau_\gamma \subseteq \tau$. Complements of γ-open sets are called γ-closed. The $\tau_\gamma$-interior [G. Sai Sundara Krishnan, 2003] of $A$ is denoted by $\tau_\gamma$-$Int(A)$ and defined to be the union of all γ-open sets of $X$ contained in $A$. The $\tau_\gamma$-closure [H. Ogata, 1991] of $A$ is denoted by $\tau_\gamma$-$Cl(A)$ and defined to be the intersection of all γ-closed sets containing $A$. A topological space $(X, \tau)$ with an operation γ on τ is said to be γ-regular [H. Ogata, 1991] if for each $x \in X$ and for each open neighborhood $V$ of $x$, there exists an open neighborhood $U$ of $x$ such that γ(U) contained in $V$. It is also to be noted that $\tau = \tau_\gamma$ if and only if $X$ is a γ-regular space [H. Ogata, 1991].

An ideal is defined as a nonempty collection $I$ of subsets $X$ satisfying the following two conditions:
1. If $A \in I$ and $B \subseteq A$, then $B \in I$.
2. If $A \in I$ and $B \in I$, then $A \cup B \in I$.

For an ideal $I$ on $(X, \tau)$, $(X, \tau, I)$ is called an ideal topological space or simply an ideal space. Given a topological space $(X, \tau)$ with an ideal $I$ on $X$ and if $P(X)$ is the set of all subsets of $X$, a set operator $(.)^* : P(X) \to P(X)$ called a local function [E. Hayashi, 1964], [K. Kuratowski, 1966] of $A$ with respect to τ and $I$ is defined as follows for a subset $A$ of $X$, $A^*(I, \tau) = \{x \in X : U \cap A \notin I$ for each neighborhood $U$ of $x\}$. A Kuratowski closure operator $Cl^*(.)$ for a topology $\tau^*(I, \tau)$, called the ∗-topology, finer than τ, is defined by $Cl^*(A) = A \cup A^*(I, \tau)$ [D. Jankovic and T. R. Hamlett, 1990]. We will simply write $A^*$ for $A^*(I, \tau)$ and $\tau^*$ for $\tau^*(I, \tau)$.

Recall that $A \subseteq (X, \tau, I)$ is called ∗-dense-in-itself [E. Hayashi, 1964] (resp. $\tau^*$-closed [D. Jankovic and T. R. Hamlett, 1990] and ∗-perfect [E. Hayashi, 1964]) if $A \subseteq A^*$ (resp. $A^* \subseteq A$ and $A = A^*$).

**Definition 2.1.** A subset $A$ of an ideal topological space $(X, \tau, I)$ is said to be
1. preopen [A. S. Mashhour, M. E. Abd El-Monsef and S. N. El-Deeb, 1982] if $A \subseteq Int(Cl(A))$.
2. pre-γ-open [H. Z. Ibrahim, 2012] if $A \subseteq \tau_\gamma$-$Int(Cl(A))$.
3. γ-preopen [G. S. S. Krishnan and K. Balachandran, 2006] if $A \subseteq \tau_\gamma$-$Int(\tau_\gamma$-$Cl(A))$.
4. γ-p-open [A. B. Khalaf and H. Z. Ibrahim, 2011] if $A \subseteq Int(\tau_\gamma$-$Cl(A))$.
5. *I*-open [D. Jankovic and T. R. Hamlett, 1992] if $A \subseteq Int(A^*)$.

6. *R-I*-open [S. Yuksel, A. Acikgoz and T. Noiri, 2005] if $A = Int(Cl^*(A))$.

7. pre-*I*-open [J. Dontchev, 1999] if $A \subseteq Int(Cl^*(A))$.

8. semi-*I*-open [E. Hatir and T. Noiri, 2002] if $A \subseteq Cl^*(Int(A))$.

9. *α-I*-open [E. Hatir and T. Noiri, 2002] if $A \subseteq Int(Cl^*(In(A)))$.

10. b-*I*-open [A. C. Guler and G. Aslim, 2005] if $A \subseteq Int(Cl^*(A)) \cup Cl^*(Int(A))$.

11. Weakly *I*-local closed [A. Keskin, T. Noiri and S. Yuksel, 2004] if $A = U \cap K$, where $U$ is an open set and $K$ is a $*$-closedset in $X$.

12. Locally closed [N. Bourbaki, 1966] if $A = U \cap K$, where $U$ is an open set and $K$ is a closed set in $X$.

**Definition 2.2.** [S. Yuksel, A. Acikgoz and T. Noiri, 2005] A point $x$ in an ideal space $(X, \tau, I)$ is called a $\delta_I$-cluster point of $A$ if $Int(Cl^*(U)) \cap A \neq \phi$ for each neighborhood $U$ of $x$. The set of all $\delta_I$-cluster points of $A$ is called the $\delta_I$-closure of $A$ and will be denoted by $\delta Cl_I(A)$. $A$ is said to be $\delta_I$-closed if $\delta Cl_I(A) = A$. The complement of a $\delta_I$-closed set is called a $\delta_I$-open set.

**Lemma 2.3.** [E. G. Yang, 2008] A subset $V$ of an ideal space $(X, \tau, I)$ is a weakly *I*-local closed set if and only if there exists $K \in \tau$ such that $V = K \cap Cl^*(V)$.

**Definition 2.4.** [E. Ekici and T. Noiri, 2009] An ideal topological space $(X, \tau, I)$ is said to be $*$-extremally disconnected if the $*$-closure of every open subset $V$ of $X$ is open.

**Theorem 2.5.** [E. Ekici and T. Noiri, 2009] For an ideal topological space $(X, \tau, I)$, the following properties are equivalent:

1. $X$ is $*$-extremally disconnected.
2. $Cl^*(Int(V)) \subseteq Int(Cl^*(V))$ for every subset $V$ of $X$.

**Lemma 2.6.** [D. Jankovic and T. R. Hamlett, 1990] Let $(X, \tau, I)$ be an ideal topological space and $A, B$ subsets of $X$. Then

1. If $A \subseteq B$, then $A^* \subseteq B^*$.
2. If $U \in \tau$, then $U \cap A^* \subseteq (U \cap A)^*$.
3. $A^*$ is closed in $(X, \tau)$.

Recall that $(X, \tau)$ is called submaximal if every dense subset of $X$ is open.

**Lemma 2.7.** [R. A. Mahmound and D. A. Rose, 1993] If $(X, \tau)$ is submaximal, then $PO(X, \tau) = \tau$.

**Corollary 2.8.** [J. Dontchev, 1999] If $(X, \tau)$ is submaximal, then for any ideal $I$ on $X$, $PIO(X) = \tau$.

Where $PIO(X)$ is the family of all pre-*I*-open subsets of $(X, \tau, I)$.

**Proposition 2.9.** [H. Ogata, 1991] Let $\gamma : \tau \to p(X)$ be a regular operation on $\tau$. If $A$ and $B$ are $\gamma$-open, then $A \cap B$ is $\gamma$-open.

## 3. Pre-*γ-I*-Open Sets

**Definition 3.1.** A subset $A$ of an ideal topological space $(X, \tau, I)$ with an operation $\gamma$ on $\tau$ is called pre-*γ-I*-open if $A \subseteq \tau_\gamma\text{-}Int(Cl^*(A))$.

We denote by $P\gamma IO(X, \tau, I)$ the family of all pre-*γ-I*-open subsets of $(X, \tau, I)$ or simply write $P\gamma IO(X, \tau)$ or $P\gamma IO(X)$ when there is no chance for confusion with the ideal.

**Theorem 3.2.** Every $\gamma$-open set is pre-*γ-I*-open.

**Proof.** Let $(X, \tau, I)$ be an ideal topological space and $A$ a $\gamma$-open set of $X$. Then $A = \tau_\gamma\text{-}Int(A) \subseteq \tau_\gamma\text{-}Int(A \cup A^*) = \tau_\gamma\text{-}Int(Cl^*(A))$.

The converse of the above theorem is not true in general as shown in the following example.

**Example 3.3.** Consider $X = \{a, b, c\}$ with $\tau = \{\phi, X, \{a, c\}\}$ and $I = \{\phi, \{b\}\}$. Define an operation $\gamma$ on $\tau$ by $\gamma(A) = X$ for all $A \in \tau$. Then $A = \{a, b\}$ is a pre-*γ-I*-open set which is not $\gamma$-open.

**Theorem 3.4.** Every pre-*γ-I*-open set is pre-*γ*-open.

**Proof.** Let $(X, \tau, I)$ be an ideal topological space and $A$ a pre-*γ-I*-open set of $X$. Then,
$A \subseteq \tau_\gamma\text{-}Int(Cl^*(A)) \subseteq \tau_\gamma\text{-}Int(Cl(A))$.

The converse of the above theorem is not true in general as shown in the following example.

**Example 3.5.** Consider $X = \{a, b, c\}$ with $\tau = \{\phi, X, \{b, c\}\}$ and $I = \{\phi, \{c\}\}$. Define an operation $\gamma$ on $\tau$ by $\gamma(A) = X$ for all $A \in \tau$. Set $A = \{c\}$, since $A^* = \phi$ and $Cl^*(A) = A$, then A is a pre-*γ*-open set which is not pre-*γ-I*-open.

**Theorem 3.6.** Every pre-*γ-I*-open set is pre-*I*-open.

**Proof.** Let $(X, \tau, I)$ be an ideal topological space and $A$ a pre-*γ-I*-open set of $X$. Then,
$A \subseteq \tau_\gamma\text{-}Int(Cl^*(A)) \subseteq Int(Cl^*(A))$.

The converse of the above theorem is not true in general as shown in the following example.

**Example 3.7.** Consider $X = \{a, b, c\}$ with $\tau = \{\phi, X, \{c\}\}$ and $I = \{\phi, \{c\}\}$. Define an operation $\gamma$ on $\tau$ by $\gamma(A) = X$ for all $A \in \tau$. Then $A = \{c\}$ is a pre-*I*-open set which is not pre-*γ-I*-open.

**Theorem 3.8.** Every pre-*γ-I*-open set is $\gamma$-preopen.

**Proof.** Let $(X, \tau, I)$ be an ideal topological space and $A$ a pre-*γ-I*-open set of $X$. Then,
$A \subseteq \tau_\gamma\text{-}Int(Cl^*(A)) \subseteq \tau_\gamma\text{-}Int(Cl(A)) \subseteq \tau_\gamma\text{-}Int(\tau_\gamma\text{-}Cl(A))$.

The converse of the above theorem is not true in general as shown in the following example.

**Example 3.9.** Consider $X = \{a, b, c\}$ with $\tau = \{\phi, X, \{b\}, \{a, b\}\}$ and $I = \{\phi, \{b\}\}$. Define an operation $\gamma$ on $\tau$ by $\gamma(A) = X$ for all $A \in \tau$. Then $A = \{b, c\}$ is a $\gamma$-preopen set which is not pre-$\gamma$-$I$-open.

**Theorem 3.10.** Every pre-$\gamma$-$I$-open set is $\gamma$-p-open.

**Proof.** Let $(X, \tau, I)$ be an ideal topological space and $A$ a pre-$\gamma$-$I$-open set of $X$. Then, $A \subseteq \tau_\gamma\text{-}Int(Cl^*(A)) \subseteq \tau_\gamma\text{-}Int(Cl(A)) \subseteq Int(\tau_\gamma\text{-}Cl(A))$.

The converse of the above theorem is not true in general as shown in the following example.

**Example 3.11.** Consider $X = \{a, b, c, d\}$ with $\tau = P(X)$ and $I = \{\phi\}$. Define an operation $\gamma$ on $\tau$ by $\gamma(A) = X$ for all $A \in \tau$. Then $A = \{c, d\}$ is a $\gamma$-p-open set which is not pre-$\gamma$-$I$-open.

**Remark 3.12.** We have the following implications but none of this implications are reversible.

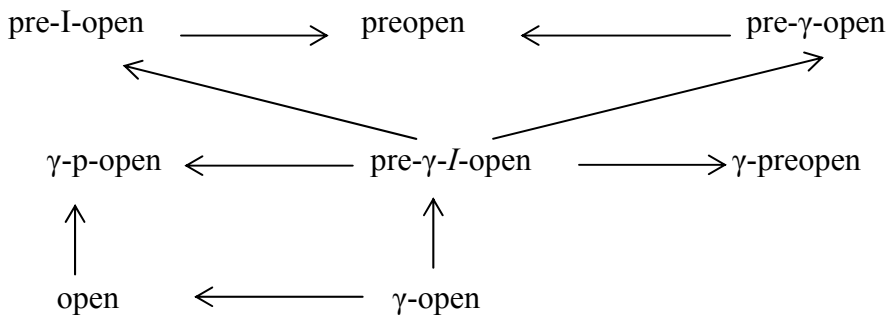

The intersection of two pre-$\gamma$-$I$-open sets need not be pre-$\gamma$-$I$-open as shown in the following example.

**Example 3.13.** Consider $X = \{a, b, c\}$ with $\tau = \{\phi, X, \{a, c\}\}$ and $I = \{\phi, \{b\}\}$. Define an operation $\gamma$ on $\tau$ by $\gamma(A) = X$ for all $A \in \tau$. Set $A = \{a, b\}$ and $B = \{b, c\}$. Since $A^* = B^* = X$, then both $A$ and $B$ are pre-$\gamma$-$I$-open. But on the other hand $A \cap B = \{b\} \notin P\gamma IO(X, \tau)$.

**Theorem 3.14.** Let $(X, \tau, I)$ be an ideal topological space and $\{A_\alpha : \alpha \in \Delta\}$ a family of subsets of $X$, where $\Delta$ is an arbitrary index set. Then,
1. If $A_\alpha \in P\gamma IO(X, \tau)$ for all $\alpha \in \Delta$, then $\cup_{\alpha \in \Delta} A_\alpha \in P\gamma IO(X, \tau)$.
2. If $A \in P\gamma IO(X, \tau)$ and $U \in \tau_\gamma$, then $A \cap U \in P\gamma IO(X, \tau)$. Where $\gamma$ is regular operation on $\tau$.

**Proof.**
1. Since $\{A_\alpha : \alpha \in \Delta\} \subseteq P\gamma IO(X, \tau)$, then $A_\alpha \subseteq \tau_\gamma\text{-}Int(Cl^*(A_\alpha))$ for each $\alpha \in \Delta$. Then we have
$\cup_{\alpha \in \Delta} A_\alpha \subseteq \cup_{\alpha \in \Delta} \tau_\gamma\text{-}Int(Cl^*(A_\alpha)) \subseteq \tau_\gamma\text{-}Int(\cup_{\alpha \in \Delta} Cl^*(A_\alpha)) \subseteq \tau_\gamma\text{-}Int(Cl^*(\cup_{\alpha \in \Delta} A_\alpha))$.
This shows that $\cup_{\alpha \in \Delta} A_\alpha \in P\gamma IO(X, \tau)$.
2. By the assumption, $A \subseteq \tau_\gamma\text{-}Int(Cl^*(A))$ and $U = \tau_\gamma\text{-}Int(U)$. Thus using Lemma 2.6, we have $A \cap U \subseteq \tau_\gamma\text{-}Int(Cl^*(A)) \cap \tau_\gamma\text{-}Int(U) = \tau_\gamma\text{-}Int(Cl^*(A) \cap U) = \tau_\gamma\text{-}Int((A^* \cup A) \cap U) = \tau_\gamma\text{-}Int((A^* \cap U) \cup (A \cap U)) \subseteq \tau_\gamma\text{-}Int((A \cap U)^* \cup (A \cap U)) = \tau_\gamma\text{-}Int(Cl^*(A \cap U))$.
This shows that $A \cap U \in P\gamma IO(X, \tau)$.

**Proposition 3.15.** For an ideal topological space $(X, \tau, I)$ with an operation $\gamma$ on $\tau$ and $A \subseteq X$ we have:
1. If $I = \{\phi\}$, then $A$ is pre-$\gamma$-$I$-open if and only if $A$ is pre-$\gamma$-open.
2. If $I = P(X)$, then $P\gamma IO(X) = \tau_\gamma$.

**Proof.** 1. By Theorem 3.4, we need to show only sufficiency. Let $I = \{\phi\}$, then $A^* = Cl(A)$ for every subset $A$ of $X$. Let $A$ be pre-$\gamma$-open, then $A \subseteq \tau_\gamma\text{-}Int(Cl(A)) = \tau_\gamma\text{-}Int(A^*) \subseteq \tau_\gamma\text{-}Int(A \cup A^*) = \tau_\gamma\text{-}Int(Cl^*(A))$ and hence $A$ is pre-$\gamma$-$I$-open. 2. Let $I = P(X)$, then $A^* = \phi$ for every subset $A$ of $X$. Let $A$ be any pre-$\gamma$-$I$-open set, then $A \subseteq \tau_\gamma\text{-}Int(Cl^*(A)) = \tau_\gamma\text{-}Int(A \cup A^*) = \tau_\gamma\text{-}Int(A \cup \phi) = \tau_\gamma\text{-}Int(A)$ and hence $A$ is $\gamma$-open. By Theorem 3.2, we obtain $P\gamma IO(X) = \tau_\gamma$.

**Remark 3.16.**
1. If a subset $A$ of a $\gamma$-regular space $(X, I, \tau)$ is open then $A$ is pre-$\gamma$-$I$-open.
2. If a subset $A$ of a submaximal space $(X, I, \tau)$ is pre-$\gamma$-$I$-open then $A$ is open.
3. If $(X, I, \tau)$ is $\gamma$-regular space and $I = P(X)$, then $A$ is pre-$\gamma$-$I$-open if and only if $A$ is open.

**Remark 3.17.** Let $(X, I, \tau)$ be a $\gamma$-regular space and $I = P(X)$. Then
1. If $A$ is $R$-$I$-open then $A$ is pre-$\gamma$-$I$-open.

2. If $A$ is $\delta_I$-open then $A$ is pre-$\gamma$-$I$-open.
3. If $A$ is regular open then $A$ is pre-$\gamma$-$I$-open.
4. If $A$ is $\delta$-open then $A$ is pre-$\gamma$-$I$-open.

**Remark 3.18.** For an ideal topological space $(X, \tau, I)$ with an operation $\gamma$ on $\tau$ and $I = P(X)$ we have:
1. If $A$ is pre-$\gamma$-$I$-open then $A$ is open.
2. If $A$ is pre-$\gamma$-$I$-open then $A$ is $\alpha$-$I$-open.
3. If $A$ is pre-$\gamma$-$I$-open then $A$ is semi-$I$-open.

**Proposition 3.19.** Let $(X, \tau, I)$ be an ideal topological space and $A$ a subset of $X$. If $A$ is closed and pre-$\gamma$-$I$-open, then $A$ is $R$-$I$-open.

**Proof.** Let $A$ be pre-$\gamma$-$I$-open, then we have $A \subseteq \tau_\gamma\text{-}Int(Cl^*(A)) \subseteq Int(Cl^*(A)) \subseteq Int(Cl(A)) \subseteq Cl(A) = A$ and hence $A$ is $R$-$I$-open.

**Remark 3.20.** Let $(X, I, \tau)$ be $\gamma$-regular space. If $A \subseteq (X, I, \tau)$ is $R$-$I$-open, then $A$ is pre-$\gamma$-$I$-open.

**Remark 3.21.** If $(X, I, \tau)$ is $\gamma$-regular space and $I = \{\phi\}$. Then
1. $A$ is pre-$\gamma$-$I$-open if and only if $A$ is preopen.
2. $A$ is pre-$\gamma$-$I$-open if and only if $A$ is $\gamma$-preopen.
3. $A$ is pre-$\gamma$-$I$-open if and only if $A$ is $\gamma$-p-open.

**Proposition 3.22.** Let $(X, \tau, I)$ be an ideal topological space and $A$ a subset of $X$. If $I = \{\phi\}$ and $A$ is pre-$\gamma$-$I$-open, then $A$ is $I$-open.

**Proof.** Let $A$ be pre-$\gamma$-$I$-open, then we have $A \subseteq \tau_\gamma\text{-}Int(Cl^*(A)) \subseteq \tau_\gamma\text{-}Int(Cl(A)) \subseteq \tau_\gamma\text{-}Int(A^*) \subseteq Int(A^*)$ and hence $A$ is $I$-open.

**Remark 3.23.** If $(X, I, \tau)$ is a $\gamma$-regular space and $A$ is $\delta_I$-open then $A$ is pre-$\gamma$-$I$-open.

**Remark 3.24.** If $(X, I, \tau)$ is $\gamma$-regular then $A$ is pre-$\gamma$-$I$-open if and only if $A$ is pre-$I$-open.

**Proposition 3.25.** If $A \subseteq (X, I, \tau)$ is $*$-perfect and pre-$\gamma$-$I$-open, then $A$ is $\gamma$-open.

**Proof.** Let $A$ be $*$-perfect, then $A = A^*$ and $A \subseteq \tau_\gamma\text{-}Int(Cl^*(A)) = \tau_\gamma\text{-}Int(A \cup A^*) = \tau_\gamma\text{-}Int(A \cup A) = \tau_\gamma\text{-}Int(A)$ and hence $A$ is $\gamma$-open.

**Remark 3.26.** If $A \subseteq (X, I, \tau)$ is $*$-perfect and pre-$\gamma$-$I$-open, then $A$ is open.

**Proposition 3.27.** If $A$ is $\tau^*$-closed in $(X, I, \tau)$ and pre-$\gamma$-$I$-open, then $A$ is $\gamma$-open.

**Proof.** Let $A$ be pre-$\gamma$-$I$-open, then $A \subseteq \tau_\gamma\text{-}Int(Cl^*(A)) = \tau_\gamma\text{-}Int(A \cup A^*) = \tau_\gamma\text{-}Int(A)$ and hence $A$ is $\gamma$-open.

**Remark 3.28.** If $A$ is $\tau^*$-closed in $(X, I, \tau)$ and pre-$\gamma$-$I$-open, then $A$ is open.

**Proposition 3.29.** If $A$ is $*$-perfect in $(X, I, \tau)$ and pre-$\gamma$-$I$-open, then $A$ is $I$-open.

**Proof.** Let $A$ be pre-$\gamma$-$I$-open, then $A \subseteq \tau_\gamma\text{-}Int(Cl^*(A)) = \tau_\gamma\text{-}Int(A \cup A^*) = \tau_\gamma\text{-}Int(A^*) \subseteq Int(A^*)$ and hence $A$ is $I$-open.

**Proposition 3.30.** If $A$ is $*$-dense-in-itself in $(X, I, \tau)$ and pre-$\gamma$-$I$-open, then $A$ is $I$-open.

**Proof.** Let $A$ be pre-$\gamma$-$I$-open, then $A \subseteq \tau_\gamma\text{-}Int(Cl^*(A)) = \tau_\gamma\text{-}Int(A \cup A^*) = \tau_\gamma\text{-}Int(A^*) \subseteq Int(A^*)$ and hence $A$ is $I$-open.

**Proposition 3.31.** If a subset $A$ of a $*$-extremally disconnected $\gamma$-regular space $(X, I, \tau)$ is $\alpha$-$I$-open then $A$ is pre-$\gamma$-$I$-open.

**Proof.** Let $A$ be $\alpha$-$I$-open, then $A \subseteq Int(Cl^*(Int(A))) \subseteq Cl^*(Int(A)) \subseteq Int(Cl^*(A)) = \tau_\gamma\text{-}Int(Cl^*(A))$ and hence $A$ is pre-$\gamma$-$I$-open.

**Proposition 3.32.** If a subset $A$ of a $*$-extremally disconnected $\gamma$-regular space $(X, I, \tau)$ is semi-$I$-open then $A$ is pre-$\gamma$-$I$-open.

**Proof.** Let $A$ be semi-$I$-open, then $A \subseteq Cl^*(Int(A)) \subseteq Int(Cl^*(A)) = \tau_\gamma\text{-}Int(Cl^*(A))$ and hence $A$ is pre-$\gamma$-$I$-open.

**Proposition 3.33.** If a subset $A$ of a $*$-extremally disconnected $\gamma$-regular space $(X, I, \tau)$ is b-$I$-open and $I = P(X)$, then $A$ is pre-$\gamma$-$I$-open.

**Proof.** Let $A$ be b-$I$-open, then $A \subseteq Int(Cl^*(A)) \cup Cl^*(Int(A)) \subseteq Int(A \cup A^*) \cup Cl^*(Int(A)) \subseteq Int(A \cup \phi) \cup Cl^*(Int(A)) \subseteq Int(A) \cup Cl^*(Int(A)) \subseteq Int(A) \cup (Int(A) \cup Int(A)^*) \subseteq Int(A) \cup Int(A)^* \subseteq Cl^*(Int(A)) \subseteq Int(Cl^*(A)) = \tau_\gamma\text{-}Int(Cl^*(A))$ and hence $A$ is pre-$\gamma$-$I$-open.

**Theorem 3.34.** Let $(X, I, \tau)$ be a $*$-extremally disconnected $\gamma$-regular ideal space and $V \subseteq X$, the following properties are equivalent:
1. $V$ is a $\gamma$-open set.
2. $V$ is $\alpha$-$I$-open and weakly $I$-local closed.
3. $V$ is pre-$\gamma$-$I$-open and weakly $I$-local closed.
4. $V$ is pre-$I$-open and weakly $I$-local closed.
5. $V$ is semi-$I$-open and weakly $I$-local closed.
6. $V$ is b-$I$-open and weakly $I$-local closed.

**Proof.** $(1) \Rightarrow (2)$: It follows from the fact that every $\gamma$-open set is open and every open set is $\alpha$-$I$-open and weakly $I$-local closed.
$(2) \Rightarrow (3)$: It follows from Proposition 3.31.
$(3) \Rightarrow (4)$, $(4) \Rightarrow (5)$ and $(5) \Rightarrow (6)$: Obvious.
$(6) \Rightarrow (1)$: Suppose that $V$ is a b-$I$-open set and a weakly $I$-local closed set in $X$. It follows that $V \subseteq Cl^*(Int(V)) \cup Int(Cl^*(V))$. Since $V$ is a weakly $I$-local closed set, then there exists an open set $G$ such that $V = G \cap Cl^*(V)$. It follows from Theorem 2.5 that $V \subseteq G \cap (Cl^*(Int(V)) \cup Int(Cl^*(V)))$
$= (G \cap Cl^*(Int(V))) \cup (G \cap Int(Cl^*(V)))$
$\subseteq (G \cap Int(Cl^*(V))) \cup (G \cap Int(Cl^*(V)))$
$= Int(G \cap Cl^*(V)) \cup Int(G \cap Cl^*(V))$
$= Int(V) \cup Int(V)$
$= Int(V)$
$= \tau_\gamma\text{-}Int(V)$.
Thus, $V \subseteq \tau_\gamma\text{-}Int(V)$ and hence $V$ is a $\gamma$-open set in $X$.

**Theorem 3.35.** Let $(X, I, \tau)$ be a $*$-extremally disconnected $\gamma$-regular ideal space and $V \subseteq X$, the following properties are equivalent:

1. $V$ is a $\gamma$-open set.
2. $V$ is $\alpha$-$I$-open and a locally closed set.
3. $V$ is pre-$\gamma$-$I$-open and a locally closed set.
4. $V$ is pre-$I$-open and a locally closed set.
5. $V$ is semi-$I$-open and a locally closed set.
6. $V$ is b-$I$-open and a locally closed set.

**Proof.** By Theorem 3.34, it follows from the fact that every open set is locally closed and every locally closed set is weakly $I$-local closed.

**Definition 3.36.** A subset $F$ of a space $(X, \tau, I)$ is said to be pre-$\gamma$-$I$-closed if its complement is pre-$\gamma$-$I$-open.

**Theorem 3.37.** A subset $A$ of a space $(X, \tau, I)$ is pre-$\gamma$-$I$-closed if and only if $\tau_\gamma\text{-}Cl(Int^*(A)) \subseteq A$.

**Proof.** Let $A$ be a pre-$\gamma$-$I$-closed set of $(X, \tau, I)$. Then $X-A$ is pre-$\gamma$-$I$-open and hence $X-A \subseteq \tau_\gamma\text{-}Int(Cl^*(X-A)) = X-\tau_\gamma\text{-}Cl(Int^*(A))$. Therefore, we have $\tau_\gamma\text{-}Cl(Int^*(A)) \subseteq A$.

Conversely, let $\tau_\gamma\text{-}Cl(Int^*(A)) \subseteq A$. Then $X-A \subseteq \tau_\gamma\text{-}Int(Cl^*(X-A))$ and hence $X-A$ is pre-$\gamma$-$I$-open. Therefore, $A$ is pre-$\gamma$-$I$-closed.

**Theorem 3.38.** If a subset $A$ of a space $(X, \tau, I)$ is pre-$\gamma$-$I$-closed, then $Cl(\tau_\gamma\text{-}Int(A)) \subseteq A$.

**Proof.** Let $A$ be any pre-$\gamma$-$I$-closed set of $(X, \tau, I)$. Since $\tau^*(I)$ is finer than $\tau$ and $\tau$ is finer than $\tau_\gamma$, we have $Cl(\tau_\gamma\text{-}Int(A)) \subseteq \tau_\gamma\text{-}Cl(\tau_\gamma\text{-}Int(A)) \subseteq \tau_\gamma\text{-}Cl(Int(A)) \subseteq \tau_\gamma\text{-}Cl(Int^*(A))$. Therefore, by Theorem 3.37, we obtain $Cl(\tau_\gamma\text{-}Int(A)) \subseteq A$.

لسەر کومێن pre-γ-I فەکرى ل ڤالاهیێن غوونەیی توبولوجی دا

**کورتی:**

ژ فێ فەکولینێ جورەکیێ ژ کوما ئەم بدەنە نیاسین وخواندن بناڤێ کومێن فەکرى ژ جورىٛ pre-γ-I ل ڤالاهیێننموونەیی توبولوجی دا.

**الملـخص:**

الغرض من هذا العمل هو تقديم و دراسة صنف من المجموعات والتي اسميناها بالمجموعات المفتوحة من النمط pre-γ-I في الفضاء التوبولوجي المثالي.